\documentclass[11pt]{article}
\usepackage{mathrsfs}
\usepackage{amscd}
\usepackage{amsmath,amsfonts,amssymb,amscd}
\usepackage{indentfirst,graphics,epsfig,psfrag}
\input{epsf}
\usepackage{ifpdf}
\usepackage{enumerate}
\usepackage{appendix}
\usepackage{enumerate}

\usepackage{lineno}

\makeatletter \@addtoreset{figure}{section} \makeatother
\makeatletter
\long\def\@makecaption#1#2{%
   \vskip 10\p@
   \setbox\@tempboxa\hbox{{#1}\ \ #2}%
   \ifdim \wd\@tempboxa >\hsize

       {#1}\ \ #2\par
   \else
       \hbox to\hsize{\hfil\box\@tempboxa\hfil}%
   \fi}
\makeatother

\newtheorem{thm}{Theorem}[section]
\newtheorem{cor}[thm]{Corollary}
\newtheorem{lem}[thm]{Lemma}

\newtheorem{pro}[thm]{Proposition}

\newcommand{\qed}{{\hfill\rule{3pt}{7pt}}}

\def\qed{\hfill \rule{4pt}{7pt}}

\begin{document}
\title{The strong rainbow vertex-connection of graphs
\footnote{Supported by NSFC.}}
\author{
\small  Xueliang Li, Yaping Mao, Yongtang Shi\\
\small Center for Combinatorics and LPMC-TJKLC\\
\small Nankai University, Tianjin 300071, P.R. China\\
\small E-mails: lxl@nankai.edu.cn; maoyaping@ymail.com; shi@nankai.edu.cn
 }
\date{}
\maketitle
\begin{abstract}
A vertex-colored graph $G$ is said to be {\it rainbow vertex-connected}
if every two vertices of $G$ are connected by a path whose internal vertices have
distinct colors, such a path is called a rainbow path. The {\it rainbow vertex-connection number} of a
connected graph $G$, denoted by $rvc(G)$, is the smallest number of
colors that are needed in order to make $G$ rainbow
vertex-connected. If for every pair $u,v$ of distinct vertices, $G$
contains a rainbow $u-v$ geodesic, then $G$ is {\it strong rainbow
vertex-connected}. The minimum number $k$ for which there exists a
$k$-vertex-coloring of $G$ that results in a strongly rainbow
vertex-connected graph is called the {\it strong rainbow
vertex-connection number} of $G$, denoted by $srvc(G)$. Observe that
$rvc(G)\leq srvc(G)$ for any nontrivial connected graph $G$.

In this paper, sharp upper and lower bounds of $srvc(G)$ are given for a
connected graph $G$ of order $n$,\ that is, $0\leq srvc(G)\leq
n-2$. Graphs of order $n$ such that $srvc(G)=1,\, 2,\, n-2$ are characterized, respectively.
It is also shown that, for each pair $a$, $b$ of integers with $a\geq 5$ and $b\geq (7a-8)/5$, there
exists a connected graph $G$ such that $rvc(G)=a$ and $srvc(G)=b$. \\[2mm]
{\bf Keywords:} vertex-coloring; rainbow vertex-connection; (strong) rainbow vertex-connection number.\\[2mm]
{\bf AMS subject classification 2010:} 05C15, 05C40, 05C76.
\end{abstract}

\section{Introduction}

All graphs considered in this paper are finite, undirected and simple.
We follow the notation and terminology of Bondy and Murty \cite{Bondy},
unless otherwise stated. Consider an edge-coloring (not necessarily
proper) of a graph $G=(V,E)$. We say that a path of $G$ is rainbow, if no
two edges on the path have the same color. An edge-colored graph $G$ is \emph{rainbow
connected} if every two vertices are connected by a rainbow path. An
edge-coloring is a \emph{strong rainbow coloring} if between every
pair of vertices, one of their geodesics, i.e., shortest paths, is a
rainbow path. The minimum number of colors required to rainbow color
a graph $G$ is called \emph{the rainbow connection number}, denoted
by $rc(G)$. Similarly, the minimum number of colors required to
strongly rainbow color a graph $G$ is called the \emph{strong
rainbow connection number}, denoted by $src(G)$. Observe that
$rc(G)\leq src(G)$ for every nontrivial connected graph $G$. The
notions of rainbow coloring and strong rainbow coloring were
introduced by Chartrand et al. \cite{Chartrand}. There are many
results on this topic, we refer to \cite{Y. Caro, CFMY}.

In \cite{M.Krivelevich}, Krivelevich and Yuster proposed a similar
concept, the concept of rainbow vertex-connection. A vertex-colored graph $G$ is
\emph{rainbow vertex-connected} if every two vertices are connected by
a path whose internal vertices have distinct colors, and such a path is called a rainbow
path. The \emph{rainbow vertex-connection number} of a connected graph $G$,
denoted by $rvc(G)$, is the smallest number of colors that are
needed in order to make $G$ rainbow vertex-connected. Note the
trivial fact that $rvc(G)=0$ if and only if $G$ is a complete graph (here
an uncolored graph is also viewed as a colored one with 0 color).
Also, clearly, $rvc(G)\geq diam(G)-1$ with equality if the diameter
is $1$ or $2$. In \cite{CLS}, the authors considered the complexity
of determining the rainbow vertex-connection of a graph. In
\cite{M.Krivelevich} and \cite{LS1}, the authors gave upper bounds
for $rvc(G)$ in terms of the minimum degree of $G$.

For more results on the rainbow connection and rainbow
vertex-connection, we refer to the survey \cite{LiSun} and a new book
\cite{LiSun1} of Li and Sun.

A natural idea is to introduce the concept of strong rainbow
vertex-connection. A vertex-colored graph $G$ is \emph{strongly rainbow
vertex-connected}, if for every pair $u,v$ of distinct vertices, there
exists a rainbow $u-v$ geodesic. The minimum number $k$ for which
there exists a $k$-coloring of $G$ that results in a strongly
rainbow vertex-connected graph is called \emph{the strong rainbow
vertex-connection number} of $G$, denoted by $srvc(G)$. Similarly, we have
$rvc(G)\leq srvc(G)$ for every nontrivial connected graph $G$.
Furthermore, for a nontrivial connected graph $G$, we have
$$
diam(G)-1\leq rvc(G)\leq srvc(G),
$$
where $diam(G)$ denotes the diameter of $G$. The following results on $srvc(G)$
are immediate from definition.

\begin{pro}\label{pro1}
Let $G$ be a nontrivial connected graph of order $n$. Then

$(a)$ $srvc(G)=0$ if and only if $G$ is a complete graph;

$(b)$ $srvc(G)=1$ if and only if $diam(G)=2$.
\end{pro}

Then, it is easy to see the following results.

\begin{cor} Let $K_{s,t}$, $K_{n_1,n_2,\ldots,n_k}$, $W_{n}$ and $P_n$
denote the complete bipartite graph, complete multipartite graph, wheel and path,
respectively. Then

$(1)$ For integers $s$ and $t$ with $s\geq 2,t \geq 1$, $srvc(K_{s,t})=1$.

$(2)$ For $k\geq 3$, $srvc(K_{n_1,n_2,\ldots,n_k})=1$.

$(3)$ For $n\geq 3$, $srvc(W_{n})=1$.

$(4)$ For $n\geq 3$, $srvc(P_n)=n-2$.

\end{cor}

\begin{figure}[!hbpt]
\begin{center}
\includegraphics[scale=0.8]{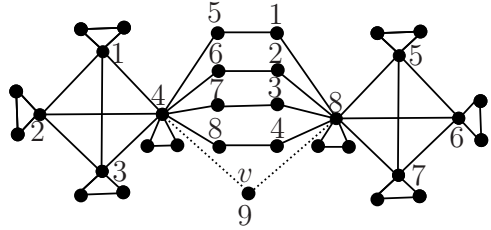}
\end{center}
\begin{center}
\caption{A counterexample for the monotonicity of the strong rainbow vertex-connection number.}
\end{center}\label{fig1}
\end{figure}

It is easy to see that if $H$ is a connected spanning subgraph of a
nontrivial (connected) graph $G$, then $rvc(G)\leq rvc(H)$. However,
the strong rainbow vertex-connection number does not have the
monotone property. An example is given in Figure \ref{fig1}, where $H=G\setminus v$ is a subgraph of $G$,
but it is easy to check that $srvc(G)=9>8=srvc(H)$. Here, one has to notice that any two cut vertices must receive
distinct colors in a rainbow vertex-coloring, just like a rainbow coloring for which any two cut edges must receive
distinct colors.

In this paper, sharp upper and lower bounds of $srvc(G)$ are given for a
connected graph $G$ of order $n$, that is, $0\leq srvc(G)\leq
n-2$. Graphs of order $n$ such that $srvc(G)=1,\, 2,\, n-2$ are characterized,
respectively. It is also shown that, for each pair $a$, $b$ of
integers with $a\geq 5$ and $b\geq (7a-8)/5$, there
exists a connected graph $G$ such that $rvc(G)=a$ and $srvc(G)=b$.

\section{Bounds and characterization of extremal graphs}

In this section, we give sharp upper and lower bounds of the strong
rainbow vertex-connection number of a graph $G$ of order $n$, that is, $0\leq
srvc(G)\leq n-2$. Furthermore, from the these bounds, we can
characterize all the graphs with $srvc(G)=0,1,n-2$, respectively.
Now we state a useful lemma.

\begin{lem}\label{lemm1}
Let $u,v$ be two vertices of a connected graph $G$. If the distance $d_G(u,v)\geq
diam(G)-1$, then there exists no geodesic containing both of $u$ and
$v$ as its internal vertices.
\end{lem}
\begin{pf}
Assume, to the contrary, that there exists a geodesic $R:w_1-w_2$
containing both $u$ and $v$ as its internal vertices. Then
$d_G(w_1,w_2)\geq d_G(u,v)+d_G(w_1,u)+d_G(v,w_2)\geq diam(G)+1$,
which contradicts to the definition of diameter.\qed
\end{pf}

\begin{thm}\label{thm1}
Let $G$ be a connected graph of order $n\ (n\geq 3)$. Then $0\leq
srvc(G)\leq n-2$. Moreover, the bounds are sharp.
\end{thm}
\begin{pf}
For $n=3$, we know $G=K_3$ or $P_3$. Since $srvc(K_3)=0<n-2$ and
$srvc(P_3)=1=n-2$, the result holds. Assume $n\geq 4$. If
$diam(G)=1$, then $G$ is a complete graph and $srvc(G)=0\leq n-2$.
If $diam(G)=2$, then we have $srvc(G)=1\leq n-2$ by Proposition 1.1.

Now suppose that $diam(G)\geq 3$. Let $diam(G)=k$ and $u,v$ be two
vertices at distance $k$. Let $P:u(=x_0),x_1,x_2,\ldots,$ $x_k(=v)$
be a geodesic connecting $u$ and $v$. Let $c$ be a
$(n-2)$-vertex-coloring of $G$ defined as: $c(u)=c(x_{k-1})=1$,
$c(x_1)=c(v)=2$, and assigning the $n-4$ distinct colors
$\{3,4,\ldots,n-2\}$ to the remaining $n-4$ vertices of $G$. Then we
will show that the coloring $c$ is indeed a strong rainbow
$(n-2)$-vertex-coloring.

It is easy to see that $d_G(u,x_{k-1})\geq k-1$ and
$d_G(v,x_{1})\geq k-1$. From Lemma \ref{lemm1}, there exists no
geodesic containing both of $u$ and $x_{k-1}$ as its internal
vertices. The same is true for vertices $v$ and $x_{1}$. So, any
geodesic connecting any two vertices of $G$ must be rainbow.
Thus, we have $0\leq srvc(G)\leq n-2$.

We show that the bounds are sharp. The complete graph $K_n$ attains
the lower bound and the path graph $P_n$ attains the upper
bound.\qed
\end{pf}\\

From Theorem \ref{thm1}, we know that $P_n$ is the graph satisfying
that $srvc(P_n)= n-2$. Actually, $P_n$ is the unique graph with this
property. That is the following theorem, which can be easily deduced
from Lemma \ref{lemm2}.

\begin{thm}
$(1)$ $srvc(G)=0$ if and only if $G$ is a complete graph;

$(2)$ $srvc(G)=1$ if and only if $diam(G)=2$;

$(3)$ $srvc(G)=n-2$ if and only if $G$ is a path of order $n$. \qed
\end{thm}

\begin{lem}\label{lemm2}
Let $G$ be a connected graph of order $n\ (n\geq 3)$. If $G$ is not
a path, then $srvc(G)\leq n-3$.
\end{lem}

\begin{pf}
For $n=3$, we know that $G=K_3$ and $srvc(G)=0=n-3$. For $n\geq 4$,
we distinguish the following two cases according to the minimum
degree $\delta(G)$ of $G$. Let $diam(G)=k$.

{\bf Case 1.}~~$\delta(G)\geq 2$.

For $4\leq n\leq 5$, $srvc(G)\leq 1\leq n-3$ since $k\leq 2$. Assume
$n\geq 6$. If $k=1$, then $G$ is a complete graph and $srvc(G)=0\leq
n-3$. If $k=2$, then it follows that $srvc(G)=1\leq n-3$.

Now suppose that $k\geq 3$ and let $P:u(=x_0),x_1,x_2,\ldots,$
$x_k(=v)$ be a geodesic of order $k$. Since $\delta(G)\geq 2$, there
exist two vertices $u'(\neq x_1)$ and $v'\neq x_{k-1}$ such that $u'$ and $v'$ are
adjacent to $u$ and $v$, respectively.

We check whether there exists a geodesic in $G$ containing both of
$u'$ and $x_{k-1}$ as its internal vertices or containing both of
$x_1$ and $v'$ as its internal vertices. If $G$ has such geodesics,
we choose one, say $Q:=w_1$-$w_2$ containing both of $u'$ and
$x_{k-1}$ as its internal vertices. It is easy to see that $w_1$ and
$w_2$ must be adjacent to $u'$ and $x_{k-1}$, respectively. We have
the following four subcases to consider.

\begin{figure}[!hbpt]
\begin{center}
\includegraphics[scale=0.8]{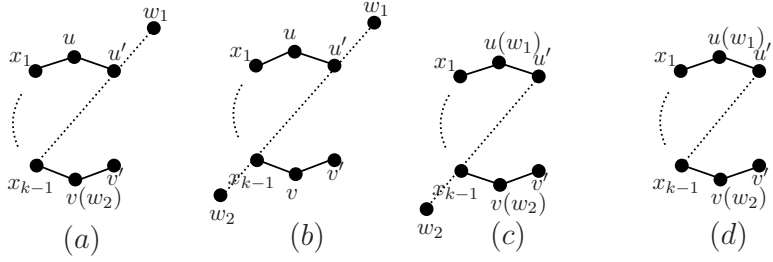}
\end{center}
\begin{center}
\caption{Four cases of $w_1$ and $w_2$}.
\end{center}
\end{figure}

{\bf Subcase 1.1.}~~$w_1\neq u$ and $w_2=v$.

Since $P$ is a geodesic, $d_G(x_1,v)\geq k-1$. Since
$d_G(u,v')+d_G(v',v)\geq d_G(u,v)=k$, we have $d_G(u,v')\geq k-1$.
By the same reason, $d_G(u',x_{k-1})\geq k-2$. Thus
$d_G(w_1,x_{k-1})\geq k-1$. By Lemma \ref{lemm1}, there exists no
geodesic containing both of $w_1$ and $x_{k-1}$ as its internal
vertices. The same is true for $v,x_{1}$ or $u,v'$. Then the
following $(n-3)$-vertex-coloring $c_1$ of $G$ is a strong rainbow
vertex-coloring: $c_1(u)=c_1(v')=1$, $c_1(v)=c_1(x_1)=2$,
$c_1(w_1)=c_1(x_{k-1})=3$ and assigning the $n-6$ distinct colors
$\{4,5,\ldots,n-3\}$ to the remaining $n-6$ vertices. So
$srvc(G)\leq n-3$.

{\bf Subcase 1.2.}~~$w_1\neq u$ and $w_2\neq v$.

It is obvious that $d_G(w_1,x_{k-1})\geq k-1$ and $d_G(w_2,u')\geq
k-1$ and $d_G(u,v)>k-1$. Then the following $(n-3)$-vertex-coloring
$c_2$ of $G$ is a strong rainbow vertex-coloring:
$c_2(u')=c_2(w_2)=1$, $c_2(w_1)=c_2(x_{k-1})=2$, $c_2(u)=c_2(v)=3$
and assigning the $n-6$ distinct colors $\{4,5,\ldots,n-3\}$ to
the remaining $n-6$ vertices. So $srvc(G)\leq n-3$.

{\bf Subcase 1.3.}~~$w_1=u$ and $w_2\neq v$.

It is easy to see that $d_G(u,x_{k-1})\geq k-1$ and
$d_G(w_2,x_1)\geq k-1$ and $d_G(u',v)\geq k-1$. From Lemma
\ref{lemm1}, we know that the following $(n-3)$-vertex-coloring
$c_3$ of $G$ is a strong rainbow vertex-coloring:
$c_3(u')=c_3(v)=1$, $c_3(u)=c_3(x_{k-1})=2$, $c_3(x_1)=c_3(w_2)=3$
and assigning the $n-6$ distinct colors $\{4,5,\ldots,n-3\}$ to
the remaining $n-6$ vertices. Hence, we have $srvc(G)\leq n-3$.

{\bf Subcase 1.4.}~~$w_1=u$ and $w_2=v$.

We will show that the following $(n-3)$-vertex-coloring $c_4$ of $G$
is a strong rainbow vertex-coloring: $c_4(u)=c_4(v)=1$,
$c_4(u')=c_4(x_{k-1})=2$, $c_4(x_1)=c_4(v')=3$ and assigning the
$n-6$ distinct colors $\{4,5,\ldots,n-3\}$ to the remaining
$n-6$ vertices.

In this case, we can use the geodesic
$P:w_1(=u),x_1,x_2,\ldots,x_{k-1},v(=w_2)$ instead of geodesic $Q$
to connect $w_1$ and $w_2$, which implies that $u'$ and $x_{k-1}$
can be assigned with the same color. From this together with
$d_G(u,v)=k$, we know that if there exists no geodesic containing
$x_1$ and $v'$ as its internal vertices, then $c_4$ is a strong
rainbow vertex-coloring.

If there exists a geodesic $R:s_1-s_2$ containing both of $x_1$ and
$v'$ as its internal vertices and $s_1\neq u$, $s_2\neq v$ or
$s_1\neq u$, $s_2=v$ or $s_1=u$, $s_2\neq v$, we can employ similar
discussions as the above three subcases of Case 1 to get $srvc(G)\leq n-3$.

For the remaining case that $s_1=u$ and $s_2=v$, we can use the geodesic
$P$ instead of geodesic $R$ to connect $s_1$ and $s_2$, which
implies that $v'$ and $x_{1}$ can be assigned with the same color.
Thus $c_4$ is indeed a strong rainbow $(n-3)$-vertex-coloring of $G$, which
results in $srvc(G)\leq n-3$.

If there is no geodesic containing both of $u'$ and $x_{k-1}$ as its internal
vertices and containing both of $v'$ and $x_{1}$ as its internal vertices,
then it is obvious that $c_4$ is also a strong rainbow vertex-coloring of $G$.
Therefore, $srvc(G)\leq n-3$.

{\bf Case 2.}~~$G$ has pendant vertices.

In this case, we will show that $srvc(G)\leq n-3$ by induction on
$n$. If $n=4$, then $G$ must be the star $K_{1,3}$ or a graph obtained
by identifying a vertex of $K_3$ to a vertex of $K_2$. From
Proposition \ref{pro1}, we have $srvc(G)=1=n-3$ since $diam(G)=2$.
Suppose that the assertion holds for a graph $G$ of smaller order. We can always find
a pendant vertex $v$ in $G$ such that $H=G-v$ is not a path. Let $u$
be adjacent to $v$ in $G$. We distinguish the following two
subcases.

{\bf Subcase 2.1.}~~$\delta(H)=1$.

Since $H$ has pendant vertices but $H$ is not a path, by induction
hypothesis, $srvc(H)\leq n-4$. We give $G$ a strong rainbow
$(n-4)$-vertex-coloring. Without loss of generality, suppose that color 1
was assigned to $u$ in $H$. Then we give $u$ a fresh color instead of $1$ and color
$v$ with $1$. Such a $(n-3)$-vertex-coloring of $G$ is a strong
rainbow vertex-coloring. Thus, we have $srvc(G)\leq n-3$.

{\bf Subcase 2.2.}~~$\delta(H)\geq 2$.

In this subcase, we can get $srvc(H)\leq n-4$ by a similar
discussion to Case 1. We also can obtain $srvc(G)\leq n-3$ by giving
a same vertex-coloring of $G$ as Subcase 2.1.

From the above arguments, we obtain that $srvc(G)\leq n-3$ if $G$ is
not a path.\qed

\end{pf}

\section{The difference of $rvc(G)$ and $srvc(G)$}

\begin{figure}[!hbpt]
\begin{center}
\includegraphics[scale=0.8]{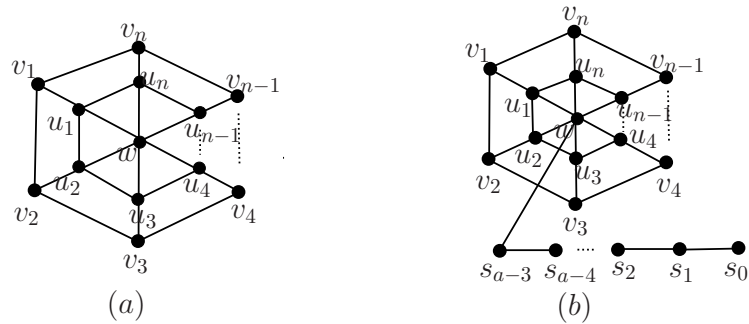}
\end{center}
\begin{center}
\caption{Graphs that are used in the proof of Theorem \ref{thm4}.}
\end{center}\label{fig3}
\end{figure}

In \cite{Chartrand}, the authors proved that for any pair $a$, $b$
of integers with $a=b$ or $3\leq a< b$ and $b\geq (5a-6)/3$, there exists a
connected graph $G$ such that $rc(G)=a$ and $src(G)=b$. Later, Chen and Li
\cite{chen-li} confirmed a conjecture of \cite{Chartrand} that for any pair $a$, $b$
of integers, there is a connected graph $G$ such that  $rc(G)=a$ and $src(G)=b$
if and only if $a=b\in \{1,2\}$ or $3\leq a\leq b$. For the two vertex-version
parameters $rvc(G)$ and $srvc(G)$, we can obtain a similar result
as follows.

\begin{thm}\label{thm4}
Let $a$ and $b$ be integers with $a\geq 5$ and $b\geq (7a-8)/5$.
Then there exists a connected graph $G$ (as shown in Figure
\ref{fig3}(b)) such that $rvc(G)=a$ and $srvc(G)=b$.\qed
\end{thm}

We construct a two-layers-wheel graph, denoted by $W^2_n$, as
follows: given two $n$-cycles $C_n^1: u_1, u_2, \ldots, u_{n}, u_1$
and $C_n^2: v_1, v_2, \ldots, v_{n}, v_1$, for $1\leq i\leq n$,
join $u_i$ to a new vertex $w$ and $v_i$ (see Figure \ref{fig3}
(a)).

\begin{lem}\label{lem3}
For $n\geq 3$, the rainbow vertex-connection number of the
two-layers-wheel $W^2_n$ is
$$rvc(W^2_n)=\left\{\begin{array}{ll}
 1 & if~n=3,\\2 & if~4\leq n\leq 6,\\3 & if~7\leq n\leq 10,\\4 & if~n\geq 11.
 \end{array}\right.$$
\end{lem}

\begin{pf}
Since $diam(W_3^2)=2$, it follows that $rvc(W_3^2)=1$. For $4\leq
n\leq 6$, $diam(W_n^2)=3$ and then $rvc(W_n^2)\geq 2$. Given a
$2$-coloring $c_1$ as follows: $c_1(w)=2$, $c_1(u_i)=1$ for $1\leq
i\leq n$, $c_1(v_i)=1$ when $i$ is odd and $c_1(v_i)=2$ otherwise.
Observe that $c_1$ is a rainbow vertex-coloring, which implies that
$rvc(W_n^2)=2$ for $4\leq n\leq 6$.

For $n=7$, $diam(W_n^2)=3$ and so $rvc(W_n^2)\geq 2$. We will show
that $rvc(W_n^2)\neq 2$. Assume, to the contrary, that
$rvc(W_n^2)=2$. Let $c'$ be a rainbow $2$-coloring of $W_n^2$. We
consider the rainbow path connecting $v_1$ and $v_4$. Since
$rvc(W_n^2)=2$, $v_2$ and $v_3$ must have distinct colors. Without
loss of generality, let $c'(v_2)=1$ and $c'(v_3)=2$. Similarly,
$v_6$ and $v_7$ must have distinct colors if considering the rainbow
path connecting $v_1$ and $v_5$. If $c'(v_6)=1$ and $c'(v_7)=2$,
then there is no rainbow path connecting $v_2$ and $v_6$ if
$c'(v_1)=2$ and also no rainbow path connecting $v_3$ and $v_7$ if
$c'(v_1)=1$. Now suppose $c'(v_6)=2$ and $c'(v_7)=1$, there is no
rainbow path connecting $v_2$ and $v_6$ if $c'(v_1)=1$. Thus
$c'(v_1)=2$. By the same reason, $v_4$ and $v_5$ must have distinct
colors. But there is no rainbow path connecting $v_4$ and $v_7$ if
$c'(v_4)=1$ and $c'(v_5)=2$ and no rainbow path connecting $v_2$ and
$v_5$ if $c'(v_4)=2$ and $c'(v_5)=1$, a contradiction. Hence, we
have $rvc(W_n^2)\geq 3$. Define the $3$-coloring $c_2$ of $W_n^2$ as
follows: $c_2(w)=3$, $c_2(u_i)=1$ when $i$ is odd and $c_2(u_i)=2$
otherwise; $c_2(v_i)=3$ when $i$ is odd and $c_2(v_i)=2$ when $i$ is
even for $1\leq i\leq 5$, $c_2(v_6)=1$, $c_2(v_7)=2$. It is easy to
check that $c_2$ is a rainbow vertex-coloring, which means that
$rvc(W_n^2)=3$ for $n=7$.

For $8\leq n\leq 9$, $diam(W_n^2)=4$ and so $rvc(W_n^2)\geq 3$. In
this case, we define the $3$-coloring $c_3$ of $W_n^2$ as follows:
$c_3(w)=3$, $c_3(u_i)=1$ when $i$ is odd and $c_3(u_i)=2$ otherwise,
$c_3(v_i)=1$ when $i\equiv 2(mod~3)$, $c_3(v_i)=2$ when $i\equiv
0(mod~3)$, $c_3(v_i)=3$ when $i\equiv 1(mod~3)$. This coloring is
also a rainbow vertex-coloring and it follows that $rvc(W_n^2)=3$
for $8\leq n\leq 9$.

For $n=10$, $diam(W_n^2)=4$ and so $rvc(W_n^2)\geq 3$. Define the
$3$-coloring $c_4$ as follows: $c_4(w)=3$, $c_4(u_i)=1$ for $1\leq
i\leq 5$ and $c_4(u_i)=2$ for $6\leq i\leq 10$, $c_4(v_1)=2$,
$c_4(v_i)=i-1$ for $2\leq i\leq 4$, $c_4(v_i)=i-4$ for $5\leq i\leq
7$, $c_4(v_8)=2$, $c_4(v_9)=1$, $c_4(v_{10})=3$. One can check that
$c_4$ is a rainbow vertex-coloring and it follows that
$rvc(W_n^2)=3$ for $n=10$.

Finally, suppose that $n\geq 11$. Observe that the $4$-coloring $c$
is a rainbow vertex-coloring: $c(w)=3$, $c(u_i)=1$ when $i$ is odd
and $c(u_i)=2$ otherwise, $c(v_i)=4$ for all $i$. It remains to show
that $rvc(W_n^2)\geq 4$. Assume, to the contrary, that
$rvc(W_n^2)=3$. Let $c'$ be a rainbow $3$-vertex-coloring of
$W_n^2$. Without loss of generality, assume that $c'(u_1)=1$. For
each $i$ with $6\leq i\leq n-4$, $v_1,u_1,w,u_i,v_i$ is the only
$v_1-v_i$ path of length $4$ in $W_n^2$ and so $u_1$, $w$ and $u_i$
must have different colors. Without loss of generality, let
$c'(w)=3$ and $c'(u_i)=2$. Since $c'(u_6)=2$, it follows that
$c'(u_n)=1$. This forces $c'(u_5)=2$, which in turn forces
$c'(u_{n-1})=1$. Similarly, $c'(u_{n-1})=1$ forces $c'(u_{4})=2$;
$c'(u_{4})=2$ forces $c'(u_{n-2})=1$; $c'(u_{n-2})=1$ forces
$c'(u_{3})=2$; $c'(u_{3})=2$ forces $c'(u_{n-3})=1$; $c'(u_{n-3})=1$
forces $c'(u_{2})=2$. There is no rainbow $v_{2}-v_{7}$ path
in $W_n^2$, which is a contradiction. Therefore, $rvc(W_n^2)=4$ for
$n\geq 11$.\qed
\end{pf}

\begin{lem}\label{lem4}
For $n\geq 3$, the strong rainbow vertex-connection number of the
two-layers-wheel $W^2_n$ is
$$srvc(W^2_n)=\left\{\begin{array}{ll}
 \lceil\frac{n}{5}\rceil& if~n=3,6;\\
 \lceil\frac{n}{5}\rceil+1& if~n\geq 4~and~n\neq 6.
 \end{array}\right.$$
\end{lem}

\begin{pf}
Since $diam(W_3^2)=2$, it follows by Proposition \ref{pro1} that
$srvc(W_3^2)=1$. If $4\leq n\leq 6$, we can check that the coloring
$c_1$ given in the proof of Lemma \ref{lem3} is a strong rainbow
$2$-vertex-coloring. So $srvc(W_n^2)\leq 2$. From this together with
$srvc(W_n^2)\geq rvc(W_n^2)=2$, it follows that $srvc(W_n^2)=2$. If
$7\leq n\leq 10$, we can check that the coloring $c_2$, $c_3$ and
$c_4$ given in the proof of Lemma \ref{lem3} is a strong rainbow
$3$-vertex-coloring. So $srvc(W_n^2)\leq 3$. Combining this with
$srvc(W_n^2)\geq rvc(W_n^2)=3$, we have $srvc(W_n^2)=3$.

Now we may assume that $n\geq 11$. Then there is an integer $k$ such
that $5k-4\leq n \leq 5k$. We first show that $srvc(W_n^2)\geq k+1$.
Assume, to the contrary, that $srvc(W_n^2)\leq k$. Let $c$ be a
strong rainbow $k$-vertex-coloring of $W_n$. If $C_n^1\cup \{w\}$
uses all the $k$ colors, it is easy to see that $w$ and $u_i$
must have distinct colors, which implies $c(u_j)\in
\{1,2,\ldots,k-1\}$ for $1\leq i\leq n$. If there exists one color
which only appears in $V(C_n^2)$, then we also have $c(u_j)\in
\{1,2,\ldots,k-1\}$ for $1\leq i\leq n$. Since $d(w)=n>5(k-1)$,
there exists one subset $S\subseteq V(C_n^1)$ such that $|S|=6$ and
all vertices in $S$ are colored the same. Thus, there exist at least
two vertices $u',\, u''\in S$ such that $d_{C_n^1}(u',u'')\geq 5$
and $d_P(u',u'')=4$, where $P:=v',u',w,u'',v''$. Since $P$ is the
only $u'-u''$ geodesic in $W_n^2$, it follows that there is no
rainbow $v'-v''$ geodesic in $W_n^2$, a contradiction.
Therefore, $srvc(W_n^2)\geq k+1$.

To show that $srvc(W_n^2)\leq k+1$, we provide a strong rainbow
$(k+1)$-vertex-coloring $c^*$: $V(W_n^2)\rightarrow \{1,2,\ldots,
k+1\}$ of $W_n^2$ defined by

$$c^*(v)=\left\{\begin{array}{ll}
 k+1 & v=w\\1, & if~v=v_i~and~i\equiv 2(mod~5),\\2 &
 if~v=v_i~and~i\equiv 3(mod~5),\\3 & if~v=v_i~and~i\equiv 4(mod~5),\\j+1 & if~v=u_i~i\in
 \{5j+1,\ldots,5j+5\}~for~0\leq j \leq k-1, \\1 & other wise.
 \end{array}\right.$$
Therefore, $srvc(W_n^2)=k+1=\lceil\frac{n}{5}\rceil+1$ for $n\geq
11$.\qed
\end{pf}\\

{\bf Proof of Theorem \ref{thm4}.} Let $n=5b-5a+10$ and let $W_n^2$
be the two-layers-wheel. Let $G$ be the graph constructed from
$W_n^2$ and the path $P_{a-1}:s_0,s_1,s_2,\ldots,s_{a-2}$ of order
$a-1$ by identifying $w$ and $s_{a-2}$ ( see Figure \ref{fig3} $(b)$).

First, we show that $rvc(G)=a$. Since $b\geq (7a-8)/5$ and $a\geq 5$,
it follows that $b>a$ and so $n=5b-5a+10>11$. By Lemma \ref{lem3},
we then have $rvc(G)=4$. Define a vertex-coloring $c$ of the graph $G$ by
$$c(v)=\left\{\begin{array}{ll}
1 & if~v=v_i~for~1\leq i\leq n,\\a & if~v=u_i~and~i~is ~odd,\\a-1 & if~v=u_i~and~i~is~even,\\i & if~v=s_i~for~1\leq i\leq a-2.
 \end{array}\right.$$
It follows that $rvc(G)\leq a$, since $c$ is a rainbow
$a$-vertex-coloring of $G$.

It remains to show that $rvc(G)\geq a$. Assume, to the contrary,
that $rvc(G)\leq a-1$. Let $c'$ be a rainbow $(a-1)$-vertex-coloring
of $G$. Since the path $s_0,s_1,s_2,\ldots,s_{a-2}(=w),u_i$ is the
only $s_0-u_i$ path in $G$, the internal vertices of this path must
be colored differently by $c'$. We may assume, without loss of
generality, that $c'(s_i)=i$ for $1\leq i\leq a-2$. For each $j$
with $1\leq j\leq 5b-5a+10$, there is a unique $s_0-v_j$ path of
length $a$ in $G$ and so $c'(u_j)=a-1$ for $1\leq j\leq 5b-5a+10$.
Consider the vertices $v_1$ and $v_{a+2}$. Since $b\geq (7a-8)/5$,
$n=5b-5a+10\geq 2a+2$ and the possible rainbow path connecting $v_1$
and $v_{a+2}$ must be $v_1,u_1,w,u_{a+2},v_{a+2}$. But it is
impossible since $c'(u_1)=c'(u_{a+2})=a-1$, which implies that there
is no $v_1-v_{a+2}$ rainbow path, contradicting our assumption that
$c'$ is a rainbow $(a-1)$-coloring of $G$. Thus, $rvc(G)\geq a$ and
then $rvc(G)=a$.

In the following, we show that $srvc(G)=b$. Since
$n=5b-5a+10=5(b-a+2)>11$, it follows from Lemma \ref{lem4} that
$srvc(W_n^2)=b-a+3$. Let $c_1$ be a strong rainbow
$(b-a+3)$-vertex-coloring of $W_n^2$. Define a vertex-coloring $c$ of the graph
$G$ by
$$c(v)=\left\{\begin{array}{ll}
 c_1(v)& if~v\in V(W_n^2);\\b-a+3+i & if~v=s_i~for~0\leq i\leq a-3.
 \end{array}\right.$$
It follows that $srvc(G)\leq b$, since $c$ is a strong rainbow
$b$-vertex-coloring of $G$.

It remains to show that $srvc(G)\geq b$. Assume, to the contrary,
that $srvc(G)\leq b-1$. Let $c^*$ be a strong rainbow
$(b-1)$-vertex-coloring of $G$. We may assume, without loss of
generality, that $c^*(s_i)=i$ for $1\leq i\leq a-2$. For each $j$
with $1\leq j\leq 5b-5a+10$, there is a unique $s_0-w-v_j$ geodesic
in $G$, implying $c^*(u_j)\in \{a-1,a,\ldots,b-1\}$. Let
$S=\{u_j:1\leq j\leq 5b-5a+10\}$. Then $|S|=5b-5a+10$ and
$|C|=b-a+1$.  Since at most five vertices in $S$ can be colored the
same, the $b-a+1$ colors in $C$ can color at most $5(b-a+1)=5b-5a+5$
vertices, producing a contradiction. Therefore, $srvc(G)\geq b$ and
so $srvc(G)=b$. \qed

\end{document}